\documentclass[11pt]{article}
\usepackage{graphicx}
\usepackage{amsmath}

\font\bigbf=cmbx10 at 16pt

\def\ds{\displaystyle}
\def\forall{\hbox{for all}~}
\def\L{{\bf L}}

\def\ve{\varepsilon}

\def\E{{\mathcal E}}

\def\R{I\!\!R}

\def\implies{\Longrightarrow}
\def\vp{\varphi}

\def\v{\vskip 1em}
\def\O{{\cal O}}

\def\dint{\int\!\!\int}

\def\C{{\cal C}}
\def\bega{\begin{array}}
\def\enda{\end{array}}
\def\begi{\begin{itemize}}
\def\endi{\end{itemize}}

\def\meas{\hbox{meas}}
\def\bel{\begin{equation}\label}
\def\eeq{\end{equation}}
\def\sqr#1#2{\vbox{\hrule height .#2pt
\hbox{\vrule width .#2pt height #1pt \kern #1pt
\vrule width .#2pt}\hrule height .#2pt }}
\def\square{\sqr74}
\def\endproof{\hphantom{MM}\hfill\llap{$\square$}\goodbreak}

\begin{document}
\title{\bigbf Representation of Dissipative Solutions to a Nonlinear 
Variational Wave Equation}
\author{Alberto Bressan and Tao Huang\\ \\ 
Department of Mathematics, Penn State University\\  University Park,
Pa. 16802, U.S.A.\\ \\ e-mails: bressan@math.psu.edu$\,$,~txh35@psu.edu}
\maketitle

\begin{abstract}
The paper introduces a new way to construct dissipative solutions to 
a second order variational wave equation. By a variable transformation, 
from the nonlinear PDE one obtains a semilinear hyperbolic 
system with sources. In contrast with the conservative case,
here the source terms are discontinuous and the discontinuities are not always crossed transversally.  Solutions to the semilinear system
are obtained by an approximation argument, relying on Kolmogorov's compactness theorem.
Reverting to the original variables, one recovers a solution to the nonlinear wave equation
where the total energy is a monotone decreasing function of time.
\end{abstract}

\section{Introduction}
\label{sec:1}
\setcounter{equation}{0}

We consider the Cauchy problem for a nonlinear wave equation
in one space dimension
\bel{1.1}
u_{tt} - c(u)\big(c(u) u_x\big)_x~=~0\,,\eeq
with initial data
\bel{1.2}
u(0,x)~=~u_0(x)\,,\qquad\qquad
u_t(0,x)~=~u_1(x)\,.\eeq
The function $c:\R\mapsto \R_+\,$,
determining the wave speeds, is assumed to be
smooth and uniformly positive.
As long as the solution remains smooth, it is well known that its energy
\bel{1.3}\E(t)~\doteq~ {1\over 2}
\int \Big[u_t^2(t,x) +c^2\big(u(t,x)\big) \,u_x^2(t,x)\Big]
\,dx\eeq
remains constant.  It is thus natural to seek global solutions
within the set of functions having bounded energy, i.e.~with $u\in H^1(\R)$ and
$u_t, u_x\in \L^2(\R)$ at a.e.~time $t$.
In this  functional space, the existence of globally defined weak solutions was 
first proved in \cite{ZZ1, ZZ2}.  For smooth initial data, 
formation of singularities was studied in \cite{GHZ}.

A different approach was introduced in  \cite{BZ},
relying on a transformation of both independent and dependent variables.
In the new variables, the equation (\ref{1.1}) is replaced by a semilinear system,
which always admits global smooth solutions (for smooth initial data).
Going back to the original variables, this method yields a family
of solutions to the original equation (\ref{1.1}), continuously depending
on the initial data (in appropriate norms).
We recall here the main properties of these global solutions, defined both forward and backward in time:
\v\begi
\item[{\bf (P)}] The solution $t\mapsto u(t,\cdot)$
takes values in $H^1(\R)$ for every $t\in\R$.
In the $t$-$x$ plane, the function $u=u(t,x)$
is  H\"older continuous
with exponent $1/2$.   The map $t\mapsto u(t,\cdot)$ is
continuously differentiable as a map with values in $\L^p_{\rm loc}$, for all
$1\leq p<2$. Moreover, it is Lipschitz continuous w.r.t.~the $\L^2$
distance, i.e.
\bel{1.4}\big\|u(t,\cdot)-u(s,\cdot)\big\|_{\L^2}~\leq~ L\,|t-s|\eeq
for some constant $L$ and all $t,s\in\R$.
The equation (\ref{1.1}) is satisfied in integral sense:
\bel{1.5}\int\!\int \Big[\phi_t\, u_t - \big(c(u) \phi\big)_x c(u)\,u_x\Big]\,dxdt ~=~0
\eeq
for all test functions
$\phi\in\C^1_c$, continuously 
differentiable with compact support in the $t$-$x$ plane.
 Concerning the initial conditions at $t=0$, the first equality in (\ref{1.2})
is satisfied pointwise, while the second holds in  $\L^p_{\rm loc}\,$
for $p\in [1,2[\,$.
\endi
\v
The solutions constructed in \cite{ZZ1, ZZ2} are {\bf dissipative}, with
 energy $t\mapsto \E(t)$ which is nonincreasing.
On the other hand, the 
solutions obtained in \cite{BZ} are {\bf conservative} in the sense that
the energy $\E(t)=\E_0$ equals a fixed constant for almost all times
$t$.
At an exceptional set of times of measure zero, one can still define
a conserved energy in terms of a positive Radon measure. However, this will not be
absolutely continuous w.r.t.~Lebesgue measure (for example, it
may contain Dirac masses).  At these particular times, the integral 
in (\ref{1.3})
accounts only for the absolutely continuous part of the energy, and is thus
strictly smaller than $\E_0$.

By putting the equation into a semilinear form, this  variable change
provides a transparent way to understand singularity formation, and construct a 
semigroup of conservative solutions continuously depending on the initial data.
A natural question, which motivated the present paper, is whether
the same variable transformation can be used to generate a semigroup of
dissipative solutions.   We recall that, in connection with the Camassa-Holm equation,
semigroups of conservative and dissipative solutions
have been constructed respectively in \cite{BC2} and \cite{BC3}, based on a similar
approach. 

To implement such a program, the main difficulty can be explained as follows.
Using characteristic variables, one obtains
a semilinear system  whose right hand side is Lipschitz continuous in the case
of conservative solutions, but  discontinuous in the case dissipative solutions.  
Because of these discontinuities,  the existence of solutions does 
not follow from general 
theory and must be studied with care.   A guiding principle is that,
if all discontinuities are crossed transversally, then the Cauchy
problem is still well posed. This is indeed what happens for the
Camassa-Holm equation \cite{BC3}. However, in connection with (\ref{1.1})  
we now encounter a ``borderline" situation, illustrated by the  system
with discontinuous right hand side
$$\bega{l} w_Y~=~
\left\{
\bega{cl}
 \cos z - \cos w&\qquad\hbox{if}~~\max\{w,z\}~<~\pi\,,\cr
 0&\qquad\hbox{if}~~\max\{w,z\}~\geq~\pi\,,\enda\right.\cr\cr
 z_X~=~
\left\{
\bega{cl}
 \cos w - \cos z &\qquad\hbox{if}~~\max\{w,z\}~<~\pi\,,\cr
 0&\qquad\hbox{if}~~\max\{w,z\}~\geq~\pi\,.\enda\right.
 \enda
$$
As $w$ approaches the discontinuity we have $w\approx \pi$, and hence
$w_Y\approx \cos z +1$.   This is strictly positive, except when
$z\approx \pm \pi$.  This lack of transversality renders the system
much harder to study.   

Our analysis shows that, assuming $c'(u)>0$ for all $u$, global dissipative solutions
of (\ref{1.1}) can indeed be constructed by solving a semilinear system
with discontinuous right hand side.
However, in contrast with \cite{BC2, BC3, BZ}, 
solutions are not obtained as the unique fixed
points of a contractive transformation.  Instead, the  existence proof 
relies here on a compactness argument, based on the Kolmogorov-Riesz
theorem \cite{HH}.  The drawback of this approach is that it does not
guarantee the uniqueness of solutions.  We speculate that the issue of uniqueness
might be resolved by the analysis of characteristics, as in \cite{BCZ, Daf}.

The paper is organized as follows.  In Section 2 we review the variable
transformation introduced in \cite{BZ} and
describe the new semilinear system with discontinuous right hand
side, which corresponds to dissipative solutions.   At the end of this section we can 
state our main results, on the global existence of solutions to the semilinear system and to the original wave equation (\ref{1.1}).
Section~3  is the core of the paper.  The discontinuous
semilinear system is here approximated by a family of Lipschitz continuous systems,
admitting unique solutions.  As the approximation parameter $\ve\to 0$,
a compactness argument yields a subsequence strongly converging to an exact solution.
In Section~4, 
returning to the original variables $u(t,x)$, we  obtain the global existence
of a weak solution to the nonlinear wave equation (\ref{1.1})-(\ref{1.2}), 
forward in time.
The proof that this solution satisfies
all properties {\rm (P)} is very similar to the one in \cite{BZ}, 
and we thus omit most of the details.

\v
\section{An equivalent semilinear system}
\setcounter{equation}{0}

We briefly review the variable transformations introduced in [BZ].
These reduce the quasilinear wave equation (\ref{1.1}) to a semilinear system,
in characteristic variables.
Throughout this section, all equations are derived assuming that the solution is smooth.
At a later stage we will prove that the same equations remain meaningful
and provide a solution to the original equation (\ref{1.1}) 
also for general initial data $(u_0, u_1)\in H^1\times \L^2$.
Define
\bel{2.1}
\left\{ \bega{rl}R&\doteq~ u_t+c(u)u_x\,,\cr
S&\doteq ~u_t-c(u)u_x\,,\enda\right.\eeq
so that
\bel{2.2}
u_t~=~{R+S\over 2}\,,\qquad\qquad u_x~=~{R-S\over 2c}\,.\eeq
If $u$ is a smooth solution of  (\ref{1.1}), these variables satisfy
\bel{2.3}\left\{
\bega{rl}
R_t-cR_x&=~{c'\over 4c}(R^2-S^2),\cr \cr S_t+cS_x&=~{c'\over 4c}(S^2-R^2).
\enda\right.
\eeq
We can thus regard $R, S$ as the densities of backward
and forward moving waves,
respectively.
Multiplying the first equation in (2.3) by $R$ and the
second one by $S$, we obtain two balance 
laws for $R^2$ and $S^2$, namely
\bel{2.4}\left\{
\bega{rl}
(R^2)_t - (cR^2)_x & =~ {c'\over 2c}(R^2S - RS^2)\, , \cr
(S^2)_t + (cS^2)_x & =~ - {c'\over 2c}(R^2S -RS^2)\,.
\enda\right.\eeq
As a consequence, for a smooth solution $u=u(t,x)$ the following quantities are conserved:
\bel{2.5}
E\doteq {1\over 2}\big(u_t^2+c^2u_x^2\big)={R^2+S^2\over 4}\,,\qquad\qquad
M\doteq -u_tu_x={S^2-R^2\over 4c}\,.\eeq
Indeed
\bel{2.6}
\left\{ \bega{rl} E_t+(c^2M)_x&=~0\,,\cr
M_t+E_x&=~0\,.\enda\right.\eeq
One can  think of $R^2/4$ and  $S^2/4$ as the energy densities associated with  
backward and forward moving waves, respectively.   Because of the sources on the right
hand sides of (\ref{2.4}), some energy
is transferred from forward to backward waves,  or viceversa.
However, the total amount of energy remains
constant.
To deal with possibly unbounded values of $R,S$,
it is convenient to introduce a new set of dependent variables:
\bel{2.7}
w\doteq 2\arctan R\,,\qquad\qquad z\doteq 2\arctan S\,,\eeq
so that
$$
R=\tan{w\over 2}\,,\qquad S=\tan{z\over 2}\,.
$$
{}From (\ref{2.3}) one derives the equations
\bel{2.8}\left\{\bega{rl}
w_t-c\,w_x&=~\ds{2\over 1+R^2}(R_t-c\,R_x)
~=~{c'\over 2c}{R^2-S^2\over 1+R^2}\,,\cr\cr
z_t+c\,z_x&=~\ds{2\over 1+S^2}(S_t+c\,S_x)
~=~{c'\over 2c}{S^2-R^2\over 1+S^2}\,. \enda\right.\eeq

\begin{figure}[ht]
\centerline{\hbox{\includegraphics[width=9cm]{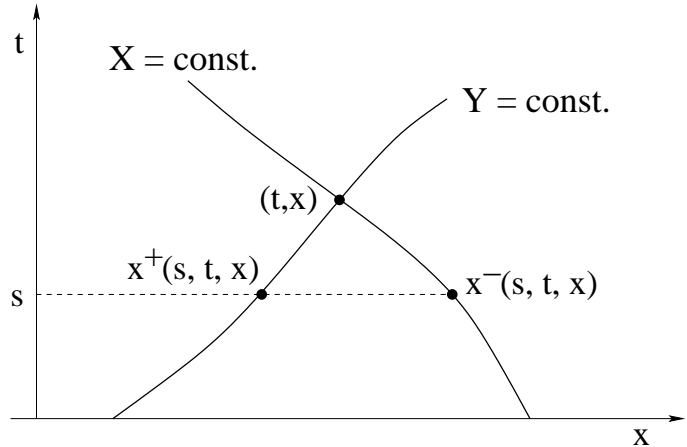}}}
\caption{\small Forward and backward characteristics through the point $(t,x)$.}
\label{f:wa17}
\end{figure}

To reduce the system to a semilinear one, we
perform a further change of independent variables.
Consider
the equations for the forward and backward characteristics  (Fig.~\ref{f:wa17}):
\bel{2.9}
\dot x^+~=~c(u)\,,\qquad\qquad \dot x^-~=~-c(u)\,.\eeq
The
characteristics passing through a given point $(t,x)$
will be denoted by
\bel{2.10}
s\mapsto x^+(s,t,x)\,,\qquad\qquad s\mapsto x^-(s,t,x)\,,\eeq
respectively.  Let an initial data $(u_0, u_1)\in H^1\times \L^2$ be given,
as in (\ref{1.2}).  In turn, this determines the functions $R, S$ in (\ref{2.1})
at time $t=0$:
\bel{RS0}\left\{\bega{rl}R(0,\cdot) &=~u_1 + c(u_0) u_{0,x}~\in \L^2(\R)\,,\cr\cr
S(0,\cdot) &=~u_1 - c(u_0) u_{0,x}~\in \L^2(\R)\,.\enda\right.\eeq
As coordinates $(X,Y)$ of a point $(t,x)$ we shall use the quantities
\bel{2.11}
X~\doteq~ \int_0^{x^-(0,t,x)} \big(1+R^2(0,x)\big)\,dx
\,,\qquad\qquad Y~\doteq~ \int_{x^+(0,t,x)}^0
\big(1+S^2(0,x)\big)\,dx\,.\eeq
Notice that this implies
\bel{2.12}
X_t- c(u)X_x ~=~0\,,\qquad\qquad Y_t+ c(u)Y_x ~=~0\,,\eeq
\bel{2.13}
(X_x)_t- (c\,X_x)_x~=~0\,,\qquad\qquad (Y_x)_t+(c\,Y_x)_x~=~0\,.\eeq
We also observe that
$$
X_x(t,x)~=~ \lim_{h\to 0}~ {1\over h} \int_{x^-(0, t, x)}
^{x^-(0, t, x+h)}\big(1+ R^2(0, x)\big)\,dx\,,
$$
$$
Y_x(t,x)~=~ \lim_{h\to 0}~ {1\over h} \int_{x^+(0, t, x)}
^{x^+(0, t, x+h)} \big(1+S^2(0, x)\big)\,dx\,.
$$
For any smooth function $f$, using  (\ref{2.12})-(\ref{2.13}) one finds
\bel{2.14}\left\{
\bega{rl}f_t+cf_x&=~f_XX_t+f_Y Y_t+cf_X X_x+cf_Y Y_x~=~ (X_t+cX_x)f_X
~=~2cX_x f_X\,,\cr
f_t-cf_x&=~f_XX_t+f_Y Y_t-cf_X X_x-cf_Y Y_x~=~ (Y_t-cY_x)f_Y
~=~-2cY_x f_Y\,.\enda\right.\eeq
We now introduce the further variables
\bel{2.15}
p~\doteq ~{1+R^2\over X_x}\,,\qquad\qquad q~\doteq ~{1+S^2\over -Y_x}\,.
\eeq
These quantities are related to the partial derivatives $X_x, Y_x$
by the identities
\bel{2.16}
(X_x)^{-1}~=~ {p\over 1+R^2}~=~p\,\cos^2{w\over 2}\,,\qquad\qquad
(-Y_x)^{-1}~=~{q\over 1+S^2}~=~q\,\cos^2{z\over 2}\,.\eeq
Notice that, if   the quantity $1+R^2$
were exactly conserved along backward characteristics, we would have
$$(1+R^2)_t -\big[c(u) (1+R^2)\big]_x~=~0\,,$$
and hence $p\equiv 1$.  In general, the variable
$p$ describes by how much the quantity $1+R^2$ fails to be
conserved along backward characteristics.  Similarly,  $q$ describes
by how much the quantity $1+S^2$ is not conserved along
forward characteristics.
\v
Starting with the nonlinear wave equation (\ref{1.1}),
using $X,Y$ as independent variables we thus obtain a semilinear
hyperbolic system with smooth coefficients for the variables $u,w,z,p,q$.
Following  \cite{BZ}, we consider the set of equations
\bel{2.17}
\left\{
\bega{rl}
w_Y&=~\theta\cdot{c'(u)\over 8c^2(u)}\,( \cos z - \cos w)\,q\,,\cr\cr
z_X&=\theta\cdot{c'(u)\over 8c^2(u)}\,( \cos w - \cos z)\,p\,,\enda\right.\eeq
\bel{2.18}
\left\{\bega{rl}p_Y&=\theta\cdot{c'(u)\over 8c^2(u)}\,
\big[\sin z-\sin w\big]\,pq\,,\cr\cr
q_X&=\theta\cdot{c'(u)\over 8c^2(u)}\,
\big[\sin w-\sin z\big]\,pq\,.\enda\right.\eeq
To obtain conservative solutions, the above equations should hold
everywhere, with $\theta\equiv 1$.
On the other hand, to construct dissipative solutions, we here choose
\bel{thdef}
\theta~=~\left\{\bega{rl} 1\quad\hbox{if}\quad \max\{w,z\}<\pi\,,\cr
0\quad\hbox{if}\quad \max\{w,z\}\geq\pi\,.\enda\right.\eeq
Finally, the function $u= u(X,Y)$ can be recovered by integrating any
of the two equations
 \bel{2.20}\left\{\bega{rl}u_Y& ={\sin z\over 4c} \, q\,,\cr
u_X &={\sin w\over 4c} \, p \,.\enda\right.\eeq
Given initial data $(u_0,u_1)\in H^1\times\L^2$ as in (\ref{1.2}),
the corresponding boundary data for the system (\ref{2.17})--(\ref{2.20})
is constructed as follows.
We first observe that the line $t=0$ corresponds
to a curve $\gamma$ in the $X$-$Y$ plane, say
$$
Y = \vp (X), \qquad X\in\R,
$$
where
$Y\doteq \vp(X)$ if and only if
$$
X=\int_0^x \big(1+R^2(0,x)\big)\,dx\,,\quad
Y=- \int_0^x \big(1+S^2(0,x)\big)\,dx\qquad\hbox{for some}~~x\in\R\,.
$$
We can use the variable $x$ as a parameter along the curve $\gamma$.
The assumptions on $u_0$ and $u_1$ imply that the corresponding 
functions
$R(0,\cdot)$ and $S(0,\cdot)$ defined at  (\ref{RS0}) are both in 
$\L^2$.   The initial energy is computed by
\bel{2.21}
\E_0~\doteq ~ {1\over 4}\int \big[R^2(0,x)+S^2(0,x)\big]\,dx 
~<~\infty .\eeq
The two functions
\bel{2.22}
X(x)\doteq \int_0^x \big(1+R^2(0,x)\big)\,dx \,,\qquad\qquad
Y(x)\doteq \int_x^0 \big(1+S^2(0,x)\big)\,dx\eeq
are well defined and absolutely continuous.  Clearly,
$X$ is strictly increasing while $Y$
is strictly decreasing. Therefore, the map $X\mapsto \vp(X)$
is continuous and strictly decreasing.  From (\ref{2.21}) it
follows
\bel{2.23}-X-4\E_0~\leq~ \vp(X)~\leq~ -X+4\E_0
\,.\eeq
As $(t,x)$ ranges over the domain $]0,\infty[\,\times\R$, the
corresponding variables $(X, Y)$ range over the  domain
\bel{2.24}
\Omega^+~\doteq~ \big\{ (X,Y)\,; ~~Y\geq \vp(X)\big\}\,.\eeq
Along the non-characteristic curve
$$
\gamma~\doteq~\big\{ (X,Y)\,;~~Y=\vp(X)\big\}\subset\R^2
$$
parameterized by $x\mapsto \big(X(x), \,Y(x)\big)$,
we can now assign the boundary data
$(\bar w, \bar z, \bar p, \bar q, \bar u)$$\in \L^\infty$
defined by
\bel{2.25}
\left\{
\bega{rcl} w(X, \vp(X))&=~
\bar w(X)   &= ~ 2\arctan R(0,x)\,,\cr
z(\vp^{-1}(Y), Y)&=~\bar z(Y) &=~  2\arctan S(0,x)\,,\cr
 u(X, \vp(X))&=~
{\bar u(X) } &=~u_0(x),
\enda\right.\qquad 
\left\{
\bega{rl}{\bar p} &\equiv  1\,,\cr
{\bar q} &\equiv  1\,.\enda\right.  \eeq
Our first main result provides the global existence of solutions to the 
discontinuous semilinear system. 
\v
{\bf Theorem 1 (existence of solutions to the
semilinear system).}  {\it Let $c=c(u)$ be a smooth function satisfying
\bel{1.6}
c(u) \geq c_0>0,\qquad\qquad c'(u)>0\qquad\qquad\forall u\in\R.\eeq
Then, for any  $(u_0, u_1)\in H^1(\R)\times \L^2(\R)$, the 
semilinear system (\ref{2.17})--(\ref{2.20}), with boundary data given by
(\ref{2.25}), (\ref{RS0}), has a solution defined for all $(X,Y)\in\Omega^+$.}
\v
In order to transform this solution back into the original variables, we set
 $f=x$ and then $f=t$ in (\ref{2.14}), and obtain 
\begin{equation}\label{th6}
\left\{
\begin{array}{ccc}
\displaystyle x_X=\frac{(1+\cos w)p}{4}\,,\cr\cr
\displaystyle  x_Y=-\frac{(1+\cos z)q}{4}\,,
\end{array}
\right.
\qquad\qquad
\left\{
\begin{array}{ccc}
\displaystyle t_X=\frac{(1+\cos w)p}{4c}\,,\cr\cr
\displaystyle  t_Y=\frac{(1+\cos z)q}{4c}\,.
\end{array}
\right.
\end{equation} 
By a direct calculation one finds $x_{XY}=x_{YX}$ and $t_{XY}=t_{YX}$. 
We can thus integrate the above equations and
recover $(t,x)$ as functions of $(X,Y)$.  
In turn, 
this yields a function $\tilde u(t,x)$ implicitly defined by
\bel{tudef}\tilde u\bigl(t(X,Y),\,x(X,Y)\bigr)~\doteq~u(X,Y) \,.\eeq

{\bf Theorem 2 (existence of dissipative solutions to the wave equation).} {\it  Let $c=c(u)$ be a smooth function satisfying (\ref{1.6}) and
consider initial data  $(u_0,u_1)\in H^1(\R)\times \L^2(\R)$. Let
$(w,z,u,p,q)$ be a solution to 
the discontinuous semilinear system (\ref{2.17})--(\ref{2.20}) with boundary data 
(\ref{2.25}).  Then the function $\tilde u(t,x)$ in (\ref{tudef}) is well defined, and 
provides a dissipative solution to the Cauchy problem (\ref{1.1})-(\ref{1.2}).
}
\v
A proof of Theorem~1 will be given in Section~3, while Theorem~2 is proved in Section~4.
\v
\section{Global solutions of the discontinuous semilinear system}
\setcounter{equation}{0}

The proof of Theorem~1 will be given in several steps.
\v
{\bf 1.} To construct solutions $(w,z,u,p,q)$ to the system (\ref{2.17})--(\ref{2.20})  we use an approximation technique.
For any $\ve>0$, consider the system of PDEs 
\bel{2.17e}
\left\{
\bega{rl}
w_Y(X,Y)&=~\theta_\ve\cdot {c'(u)\over 8c^2(u)}\,( \cos z - \cos w)\,q+\ve\,,\cr\cr
z_X(X,Y)&=~\theta_\ve\cdot {c'(u)\over 8c^2(u)}\,( \cos w - \cos z)\,p+\ve\,,\enda\right.
\eeq
\bel{2.18e}
\left\{\bega{rl}p_Y&=~\theta_\ve\cdot {c'(u)\over 8c^2(u)}\,
\big[\sin z-\sin w\big]\,pq\,,\cr\cr
q_X&=~\theta_\ve\cdot{c'(u)\over 8c^2(u)}\,
\big[\sin w-\sin z\big]\,pq\,,\enda\right.
\eeq
\bel{2.20e}u_Y~ =~{\sin z\over 4c(u)} \, q\,.\eeq
Here the coefficient $\theta_\ve = \theta_\ve (\max\{w,z\})$ is defined by setting
\bel{thetadef}\theta_\ve~\doteq~\left\{\bega{rl} 1\quad &\hbox{if}\quad \max\{z,w\}\leq \pi,\cr
 0\quad &\hbox{if}\quad \max\{z,w\}\geq \pi+\ve^3,\enda\right.\eeq
and by requiring $\theta_\ve$ to be an affine function 
of $\max\{z,w\}$ on the interval $[\pi, \pi+\ve^3]$.
Let initial data (\ref{2.25}) be given along the curve $\gamma= \bigl\{
(X,Y)\,;~ Y=\vp(X)\bigr\}
\subset\R^2$.
For convenience, we extend it to the outer region $\{(X,Y);~ Y<\vp(X)\}$
by letting $u,w,p$ be constant along vertical lines (where $X$ is constant) 
and $z,q$ be constant
along horizontal lines (where $Y$ is constant). We observe that, for any $\ve>0$,  the right hand sides of 
(\ref{2.17e})-(\ref{2.18e}) and (\ref{2.20e}) are Lipschitz continuous. Hence,
given the initial data 
\bel{ide}
\left\{\bega{rl} w(X,\vp(X))&=~\bar w(X), \cr z(Y,\vp^{-1}(Y))&=~\bar z(Y),\enda\right.
\qquad \left\{\bega{rl} p(X,\vp(X))&=~1, \cr q(Y,\vp^{-1}(Y))&=~1,\enda\right.
 \qquad u(X,\vp(X)) = \bar u(X)),
\eeq
 this semilinear hyperbolic system 
admits a unique local solution, say $(w_\ve, z_\ve, p_\ve, q_\ve, u_\ve)$. 
Indeed, this solution can be obtained as the unique fixed point of an integral transformation.  

We claim   
that this solution is globally defined on the entire
domain $\Omega^+=\{(X,Y);~Y\geq\vp(X)\}$.  This is not immediately obvious
because the equations (\ref{2.18e}) are quadratic w.r.t.~$p,q$.   
For a given $M>0$, consider the domain
\bel{OM}\Omega_M~\doteq~\{ (X,Y)\,;~Y\geq\vp(X)\,,~~X\leq M\,,~~~Y\leq M\}~\subset
~\Omega^+.\eeq
To achieve the global existence on $\Omega_M$, it suffices to  prove the uniform 
a priori bounds 
\bel{pqb}
\left\{\bega{l}
0~<~C^{-1}~\leq ~p(X,Y)~\leq ~C,\cr
0~<~C^{-1}~\leq ~q(X,Y)~\leq ~C,\enda\right.
\qquad\qquad\forall (X,Y)\in \Omega_M\,,\eeq
for some constant $C$ depending on $M$.

Observe that  (\ref{2.18e}) implies
$p_Y+q_X=0$.
Hence the differential form $pdX-qdY$ has zero integral 
along any closed curve in $\Omega_{M}$. For any $(X,Y)\in\Omega_M$, 
consider the closed curve $\Gamma=\Gamma_1\cup\Gamma_2\cup\Gamma_3$, 
where $\Gamma_1$ is the portion of boundary $\gamma$ between 
$(\vp^{-1}(Y), Y)$ and $(X, \vp(X))$,
$\Gamma_2$ is the vertical segment from $(X, \vp(X))$ to $(X, Y)$ 
and $\Gamma_3$ is the horizontal segment from $(X,Y)$ to $(\vp^{-1}(Y),Y)$. 
Integrating $pdX-qdY$ on $\Gamma$ and using the initial data (\ref{ide}), one obtains
\bel{th19}\bega{l}\ds
\int_{\vp^{-1}(Y)}^Xp(X',Y)\,dX'+\int_{\vp(X)}^Yq(X,Y')\,dY'~=~X-\vp^{-1}(Y)+Y-\vp(X)
\cr\cr
\qquad\ds \leq~ 2\left(|X|+|Y|+4\mathcal{E}_0\right),
\enda \eeq
because of  (\ref{2.23}). 
To fix the ideas, assume 
\begin{equation}\label{th3}
C_0~\doteq~\sup_u\frac{c'(u)}{8c^2(u)}~<~+\infty.
\end{equation}
Integrating the first equation of (\ref{2.18e}) along line segment from $(X, \vp(X))$ to $(X, Y)$ and using  (\ref{2.23}), (\ref{th19}), since  $p,q>0$ we obtain
\begin{eqnarray}
p(X,Y)&=&\exp\left(\int_{\vp(X)}^Y\theta\frac{c'}{8c^2}(\sin z-\sin w)q(X,Y')\,dY'\right)\notag\\
&\leq&\exp\left(2C_0\int_{\vp(X)}^Yq(X,Y')\,dY'\right)\label{th4}\\
&\leq& \exp\left( 8C_0\left(M+2\mathcal{E}_0\right)\right).\notag
\end{eqnarray}
In the same way, we also obtain
\begin{eqnarray}
p(X,Y)&=&\exp\left(-2C_0\int_{\vp(X)}^Yq(X,Y')\,dY'\right)\label{th5}\\
&\geq& \exp\left( -8C_0\left(M+2\mathcal{E}_0\right)\right).\notag
\end{eqnarray}
Together (\ref{th4}) and (\ref{th5}) yield the first estimate in (\ref{pqb}).
The second estimate is entirely similar.

In turn, this implies that 
the maps 
$$Y~\mapsto~w(X,Y),\qquad\qquad X~\mapsto~z(X,Y)$$
are uniformly Lipschitz continuous on the domain $\Omega_M$, 
say with Lipschitz constant $L$.
\v
{\bf 2.} 
Our next goal is to show that, as $\ve\to 0$, this sequence of approximations is
compact in $\L^1_{loc}(\R^2)$. 
For this purpose, some a priori estimates are needed.
Fix $\ve>0$ and consider the corresponding  solution  $(w, z, p, q, u)$
of 
(\ref{2.17e})--(\ref{thetadef}).  To shorten notation, in the following we drop the
the subscript $\ve$. 
   Define the maps
$$\bega{l} X\mapsto Y^\pi(X)~\doteq~\inf\{ Y\in [0, M]\,;~~w(X,Y) = \pi\},\cr
Y\mapsto X^\pi(Y)~\doteq~\inf\{ X\in [0, M]\,;~~z(X,Y) = \pi\},\enda$$
Observe that,   for $\ve>0$ small enough, if $Y^\pi(X)<M$, then 
\bel{wpie}w(X, Y)~\geq~ \pi+\ve^3\qquad\qquad \forall Y~\geq~ Y^\pi(X) +\ve.\eeq
Indeed,  for $Y\in [Y^\pi(X), Y^\pi(X) +\ve]$ we have
\bel{wpee}w_Y(X,Y)~\geq~\theta_\ve\cdot {c'(u)\over 8c^2(u)}( 1-\cos w) q +\ve
~\geq~{\ve\over 2}\,.\eeq
As soon as $w$ becomes greater than $  \pi + \ve^3$ we have $\theta_\ve=0$ 
and $w_Y=\ve$.
This implies (\ref{wpie}).
\v
\begin{figure}[ht]
\centerline{\hbox{\includegraphics[width=10cm]{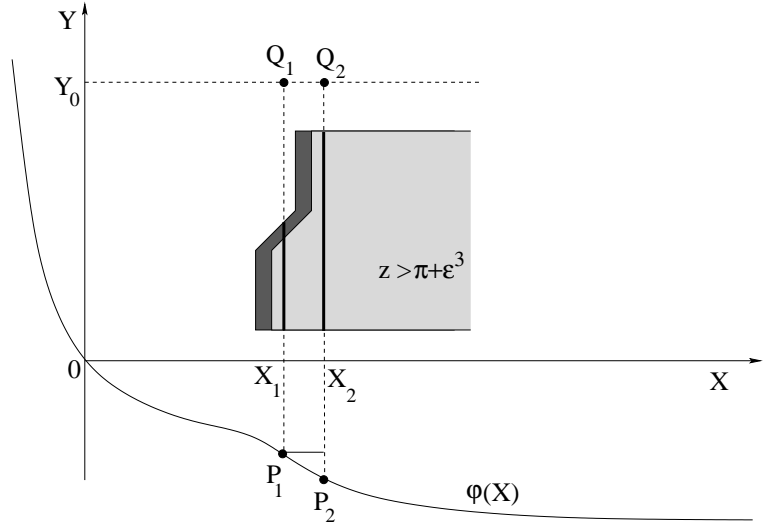}}}
\caption{\small Estimating the difference $w(X_2,Y_0)-w(X_1, Y_0)$.}
\label{f:wa18}
\end{figure}

{\bf 3.} The heart of the proof is provided by the next lemma.

{\bf Lemma 1.}
{\it  Given $M>0$ there exists a constant $C$ such that, for any $\ve\in \,]0,1]$,
the solution of (\ref{2.17e})--(\ref{ide}) satisfies the following estimates. 
For any $X_1<X_2\leq M$ and any $Y_0\in [\vp(X_1)\,~M]$, one has
\bel{kes1}\bega{l}
|w(X_1,Y_0)-w(X_2, Y_0)|+|p(X_1,Y_0)-p(X_2, Y_0)|+|u(X_1,Y_0)-u(X_2, Y_0)|\cr\cr
\quad \leq~C \bigg\{|\bar w(X_1)-\bar w(X_2)|+ |\bar u(X_1)-\bar u(X_2)|
+ |X_1-X_2|+ |\vp(X_1)-\vp(X_2)| \cr\cr
\qquad   + \meas
\Big( \{ Y\leq M\,;~~ [X_1, X_2]\cap[X^\pi(Y),\, X^\pi(Y) +\ve]~\not=~\emptyset\} \Big)\bigg\}^{1/2}.
\enda\eeq
Moreover, for $Y_1<Y_2$ and any $X_0\in [ \vp^{-1}(Y_1),\, M]$, one has
\bel{kes2}\bega{l}
|z(X_0,Y_1)-z(X_0, Y_2)|+|q(X_0,Y_1)-q(X_0, Y_2)|\cr\cr
\qquad \leq~C \,\bigg\{|\bar z(Y_1)-\bar z(Y_2)|+ |\bar u(Y_1)-\bar u(Y_2)|
+|Y_1-Y_2| + |\vp^{-1}(Y_1)-\vp^{-1}(Y_2)|\cr\cr
\qquad \qquad  +\meas
\Big( \{ X\leq M\,; ~~ [ Y_1, \, Y_2]\cap [Y^\pi(X),\, Y^\pi(X)+ \ve]~\not=~\emptyset\} \Big)\bigg\}^{1/2}.
\enda\eeq
}

The above estimates can be explained with the aid of Figure~\ref{f:wa18}.
For i=1,2, consider the points $P_i=(X_i, \vp(X_i))$, $Q_i=(X_i, Y_0)$.   
Then $w(Q_i)$ 
can be computed by solving the ODE in (\ref{2.17e})
on the interval $Y\in [\vp(X_i), Y_0]$, with initial data $w=\bar w$
at $Y=\vp(X_i)$. For $i=1,2$, the 
right hand sides of these ODEs are almost the same, except in the case where
$z(X_1,Y)<\pi$ but $z(X_2,Y)>\pi$. This motivates the presence of the last term
on the right hand side of (\ref{kes1}) and (\ref{kes2}).

{\bf Proof of Lemma 1.}  
For notational convenience, we lump together different variables and write
$$ \alpha(X,Y) ~\doteq~(w,p,u)(X,Y),\qquad\qquad  \beta(X,Y) ~\doteq~(z,q)(X,Y).$$

We first consider the easier case where
\bel{C1}w(X_i, Y)~<~ \pi\qquad\qquad \forall Y\leq Y_0,~~~i=1,2\,.\eeq

Since the maps $X\mapsto \beta(X,Y)$ are uniformly bounded and Lipschitz continuous
on bounded sets, we can find a constant $\kappa$ such that 
\begi
\item[(i)]  If 
 $[X_1, X_2]\cap[X^\pi(Y),\, X^\pi(Y) +\ve]~=~\emptyset$, then 
\bel{PY1}
\frac{\partial}{\partial Y} \Big|\alpha(X_1,Y)-\alpha(X_2,Y)\Big|~\leq ~
\kappa\Big(
|X_1-X_2|+|\alpha(X_1,Y)-\alpha(X_2,Y)|\Big).
\eeq
\item[(ii)]   If 
 $[X_1, X_2]\cap[X^\pi(Y),\, X^\pi(Y) +\ve]~\not=~\emptyset$, then 
\bel{PY2}
\frac{\partial}{\partial Y} \Big|\alpha(X_1,Y)-\alpha(X_2,Y)\Big|~\leq ~\kappa.\eeq
\endi
The differential inequalities (\ref{PY1})-(\ref{PY2}) are complemented by the estimate 
on the initial data 
\bel{ap12}\bega{l}\Big|\alpha(X_1, \vp(X_1))- \alpha(X_2, \vp(X_1))\Big|~\leq~
\Big|\bar w(X_1)- \bar w(X_2)\Big|\cr\cr
\qquad\qquad +\Big|\bar p(X_1)- \bar p(X_2)\Big|
+
\Big|\bar u(X_1)- \bar u(X_2)\Big|+
\kappa|\vp(X_1)-\vp(X_2)|\,,\enda\eeq
for a suitable constant $\kappa$. 

Using the differential inequalities  (\ref{PY1})-(\ref{PY2}) on the interval 
$Y\in [\vp(X_1), Y_0]$ together with
 (\ref{ap12}), by
a Gronwall-type estimate we obtain   
\bel{kes3}\bega{l}
|w(X_1,Y_0)-w(X_2, Y_0)|+|p(X_1,Y_0)-p(X_2, Y_0)|+|u(X_1,Y_0)-u(X_2, Y_0)|\cr\cr
\quad \leq~C_1 \bigg\{|\bar w(X_1)-\bar w(X_2)|+ |\bar u(X_1)-\bar u(X_2)|
+ |X_1-X_2|+ |\vp(X_1)-\vp(X_2)| \cr\cr
\qquad   + \meas
\Big( \{ Y\leq M\,;~~ [X_1, X_2]\cap[X^\pi(Y),\, X^\pi(Y) +\ve]~
\not=~\emptyset\} \Big)\bigg\},
\enda\eeq
for a suitable constant $C_1$.
In the case where 
\bel{C2}z(X, Y_i)~\leq~ \pi\qquad\qquad \forall X\leq X_0,~~~i=1,2\,,\eeq
an entirely similar argument yields 
\bel{kes4}\bega{l}
|z(X_0,Y_1)-z(X_0, Y_2)|+|q(X_0,Y_1)-q(X_0, Y_2)|\cr\cr
\qquad \leq~C \,\bigg\{|\bar z(Y_1)-\bar z(Y_2)|+ |\bar u(Y_1)-\bar u(Y_2)|
+|Y_1-Y_2| + |\vp^{-1}(Y_1)-\vp^{-1}(Y_2)|\cr\cr
\qquad \qquad  +\meas
\Big( \{ X\leq M\,; ~~ [ Y_1, \, Y_2]\cap [Y^\pi(X),\, Y^\pi(X)+ \ve]~\not=~\emptyset\} \Big)\bigg\}.
\enda\eeq
\v
We now study the more difficult case where (\ref{C1}) does not hold.
To fix the ideas, assume that, for some $Y_*\in [\vp(X_1), Y_0]$, we have
\bel{CS}
w(X_i, Y)~<~ \pi\qquad \forall Y< Y_*,~~~i=1,2\,,
\qquad w(X_1, Y_*) ~=~ \pi.\eeq
For $Y_0\geq Y_*$, using the triangle inequality we can write
$$\bega{rl}&|\alpha(X_1,Y_0)-\alpha(X_2,Y_0)|\cr\cr
&\quad \leq~|\alpha(X_1,Y_0)-\alpha(X_1, Y_*)|+
 |\alpha(X_1,Y_*)-\alpha(X_2,Y_*)|+ 
 |\alpha(X_2,Y_*)-\alpha(X_2,Y_0)|\cr\cr
 &\quad =~A_1+A_2+A_3.\enda$$
 Recalling (\ref{wpie}) we have the bound 
$$A_1~\leq ~  \kappa \ve $$
for some constant $\kappa$.
Moreover, by the previous arguments we already know that the estimate
(\ref{kes3}) holds when $Y_0$ is replaced by $Y_*$.   This yields a bound on $A_2$.

We now work toward an estimate of $A_3$.
Call 
$$Y^*~\doteq~\sup\{ Y\in [Y_*, Y_0]\,;~~~w(X_2, Y)<\pi\}\,.$$
Since 
$$|\alpha(X_2,Y_0)-\alpha(X_2, Y^*)|~\leq~\kappa\ve,$$ it suffices to
estimate the difference $|\alpha(X_2,Y^*)-\alpha(X_2, Y_*)|$.

The relevant equations in (\ref{2.17e})--(\ref{2.20e})
are
\bel{te}
\left\{\bega{rl}w_Y&\ds=~{c'(u)\over 8c^2(u)}(\cos z - \cos w)q+\ve\,,
\cr\cr
p_Y&=\ds~{c'(u)\over 8c^2(u)}(\sin z-\sin w)pq\,,
\cr\cr
u_Y&=\ds~{\sin z\over 4 c(u)}\,q\,,\enda\right.\eeq
with an initial data  $w(Y_*,X_2)\approx \pi$.
Since in these computations $X=X_2$ is fixed, we shall omit this variable
and write $w(Y)=w(X_2,Y)$, $z(Y)=z(X_2,Y)$, etc$\ldots$~~ 
Roughly speaking, two cases can  occur:
\begi
\item[(i)] $\sin z(Y) \approx 0$.   
In this case  $w_Y, p_Y, u_Y \approx 0$.  Hence all these functions remain almost constant.
\item[(ii)]  $\sin z(Y)$ is not close to zero.   In this case $\cos z(Y)$ is
much bigger than  $-1$, hence $w_Y(Y)$ is strictly positive.   Therefore,
$w(\cdot)$ will increase, reaching $\pi+\ve^3$ within a short time.   
After this happens, 
$|\alpha_Y(\cdot)|\leq 2\ve$, hence $\alpha$ remains almost constant.
\endi
In both cases, 
the difference $|\alpha(Y)-\alpha(Y_*)|$ remains small.

Relying on the previous ideas,  we now work out a rigorous proof. Set
$$\delta~\doteq~\pi - w(Y_*).$$
Notice that $\delta$ is bounded by the right hand side of (\ref{kes3}).
Choose constants $0<c_0<C_0$ such that
\bel{C00}c_0~\leq~{c'(u)\over 8 c(u)} ~\leq~C_0\,.\eeq
From the first equation in (\ref{te}) we deduce the lower bound
\bel{lbw}
w(Y)~\geq~w(Y_*)  - C_0 (1-\cos(\pi-2\delta))
\qquad\qquad Y\in [Y_*, Y^*]\,.\eeq
In addition, we have
\bel{e4}\left\{
\bega{rl}
|p_Y|&\leq~C_0\Big( |\sin z(Y)|+ |\sin 2\delta|\Big)+\ve\,,\cr\cr
|u_Y|&\leq~C_0\,|\sin z(Y)|\,.
\enda\right.\eeq
\bel{e5}\delta~\geq~w(Y^*)-w(Y_*)~\geq~
\int_{Y_*}^{Y^*}{c'(u)\over 8c^2(u)} (\cos z(Y)- \cos(\pi-2\delta))\, dY
.\eeq
\bel{e6}
\int_{Y_*}^{Y^*}{c'(u)\over 8c^2(u)} \Big(\cos z(Y) 
+ 1 -\delta^2
\Big)\, 
dY~\leq~\delta
.\eeq
Using (\ref{C00}),  from (\ref{e6}) we  deduce 
\bel{e7}
\bega{l}\ds
\int_{Y_*}^{Y^*}|\sin z(Y)|\, dY~=~ \int_{Y_*}^{Y^*}2\Big|\sin {z(Y)\over 2}\Big|
\,\Big|\cos {z(Y)\over 2}\Big|\, dY\cr\cr
\ds\leq~2|Y^*-Y_*|^{1/2}\left(\int_{Y_*}^{Y^*}\cos^2 {z(Y)\over 2}\, dY\right)^{1/2}
~\leq~C'\left(\int_{Y_*}^{Y^*}1+\cos z(Y)\, dY\right)^{1/2}
~\leq~ C \delta^{1/2}
\enda
\eeq
for some constants $C',C$.
By (\ref{e4}) we thus have
\bel{e8}
|p(Y^*)-p(Y_*)|~\leq~\int_{Y_*}^{Y^*}|p_Y|\, dY
~\leq~  C(\delta^{1/2}+\ve )
\eeq
and 
\bel{e9}
|u(Y^*)-u(Y_*)|~\leq~\int_{Y_*}^{Y^*}|u_Y|\, dY
~\leq~  C\delta^{1/2},
\eeq
possibly with a larger constant $C$.  Since $\delta$ is bounded by the right hand side of (\ref{kes3}), by a suitable choice of $C$ we obtain  (\ref{kes1}).      

Using (\ref{kes4}), a similar argument yields (\ref{kes2}).
\endproof 
\v
{\bf 4.} Recalling the definition  of the 
domain $\Omega_M$ at  (\ref{OM}), 
consider any rectangle $[a,b]\times [c,d]\subset \Omega_M$, and let $(\xi,\zeta)$
be any vector such that $[a+\xi,\,b+\xi]\times [c+\zeta,d+\zeta]\subset \Omega_M$.

Consider any solution of (\ref{2.17e})--(\ref{ide}), for some $\ve\in \,]0,1]$.
Since the components $w,p,u$ are uniformly 
Lipschitz continuous w.r.t.~$Y$,
we have the easy estimate
\bel{es5}\int_a^b \int_c^d|\alpha( X,Y) - \alpha(X, Y+\zeta)|\, dY\, dX~\leq~
\int_a^b C\,|\zeta|\, dX ~\leq~C(b-a)|\zeta|\,.\eeq
for some constant $C$.

Next, using (\ref{kes1}) with  $X_1\doteq X$, $X_2\doteq X+\xi$, 
we obtain
\bel{k2}\bega{l}\ds\int_a^b \int_c^d\left|\alpha (X,Y)-\alpha(X+\xi, Y)\right|\,dY\,  dX
\cr\cr
\ds\leq~C\,(d-c)\int_a^b \bigg\{\bigl|\bar w(X)-\bar w(X+\xi)\bigr|+ \bigl|\bar u(X)-\bar u(X+\xi)\bigr|
+ |\xi|+ \bigl|\vp(X)-\vp(X+\xi)\bigr| \cr\cr
\qquad   + \meas
\Big( \{ Y\leq M\,;~~ [X, X+\xi]\cap[X^\pi(Y),\, X^\pi(Y) +\ve]~\not=~\emptyset\} \Big)
\bigg\}^{1/2}\, dX\cr\cr
\ds\leq ~C'\int_a^b \Big(A(X)^{1/2} + B(X)^{1/2}\Big)\, dX\,,
\enda\eeq
where
$$A(X)~\doteq~\bigl|\bar w(X)-\bar w(X+\xi)\bigr|+ \bigl|\bar u(X)-\bar u(X+\xi)\bigr|
+ |\xi|+ \bigl|\vp(X)-\vp(X+\xi)\bigr|,$$
$$B(X)~\doteq~ \meas
\Big( \{ Y\leq M\,;~~ [X, X+\xi]\cap[X^\pi(Y),\, X^\pi(Y) +\ve]~\not=~\emptyset\} \Big)
.$$
Since the initial data $\bar w$ is bounded and measurable, while $\bar p\equiv 1$ and 
$\bar u,\vp$ are continuous, there exists some modulus of continuity
$\psi$ such that
\bel{k3}\int_a^b \Big(\bigl|\bar w(X)-\bar w(X+\xi)
\bigr|+\bigl|\bar u(X)-\bar u(X+\xi)\bigr|+|\xi|+
\bigl|\vp(X)-\vp(X+\xi)\bigr|\Big) \, dX~
\leq~\psi(|\xi|)\,.\eeq
Calling $A_0\doteq {\psi(|\xi|)\over b-a}$, we have
\bel{j1} \bega{rl}\ds\int_a^b A(X)^{1/2}\, dX&\ds=~\int_{[a,b]\cap\{ A\leq A_0\}}A(X)^{1/2}\, dX
+ \int_{[a,b]\cap\{ A>A_0\}}A(X)^{1/2}\, dX\cr\cr
&\ds\leq~\sqrt{\psi(|\xi|)\cdot( b-a)}~ + ~{1\over \sqrt{A_0}}\int_a^b A(X)\, dX~\leq~
2\sqrt{\psi(|\xi|)\cdot( b-a)}.\enda
\eeq
To estimate the integral of $B^{1/2}$ we observe that
\bel{k5}\int_a^b B(X)\, dX~\leq~(d-c)(|\xi|+\ve).\eeq
The same arguments as in (\ref{j1}), with $\psi(|\xi|)$ replaced by the right hand side of (\ref{k5}),
 now
 yield
\bel{j2} \int_a^b B(X)^{1/2}\, dX~\leq~2\sqrt{ (d-c)(|\xi|+\ve)\cdot (b-a)}.
\eeq
Together, the estimates (\ref{es5}) and (\ref{j1})-(\ref{j2}) 
yield a bound of the form
\bel{k6}\bega{l}\ds\int_a^b \int_c^d\left|\alpha (X,Y)-\alpha(X+\xi, Y+\zeta)\right|\,dY\,  dX
~\leq~\Phi(|\xi|+|\zeta|+\ve)\enda\eeq
for some continuous function $\Phi$, with $\Phi(0)=0$.
Entirely similar estimates hold for the functions $\beta=(z,q)$.
\v
{\bf 5.} Given any 
sequence $\ve_n\to 0$, consider the corresponding approximate solutions 
$U^\ve\doteq (u^\ve ,w^\ve,z^\ve,p^\ve,q^\ve)$.   
In order to use the Kolmogorov-Riesz
compactness theorem \cite{HH} and prove that a subsequence $U^{\ve_n}$ admits a subsequence
converging in $\L^1_{loc}(\Omega^+)$, the following property must be proved:
\v
{\bf (P)} {\it For any 
$\epsilon>0$ there exists $\rho>0$ such that, on any rectangle $Q\subset\Omega^+$,
one has    
\bel{Uee}\int_Q |U^{\ve_n}(X+\xi, \,Y+\zeta)- U^{\ve_n}(X,Y)|\, dX\,dY~\leq~\epsilon\eeq
whenever   $|\xi|+|\zeta|\leq \rho$, $n\geq 1$.}   

By the previous step, we have an estimate of the form
\bel{U2}\int_Q 
\bigl|U^{\ve_n}(X+\xi, \,Y+\zeta)- U^{\ve_n}(X,Y)\bigr|\, dX\,dY~\leq~
\Phi(|\xi|+|\zeta|+\ve_n).\eeq
Choose $\rho'>0$ small enough so that 
$\Phi(2\rho')<\epsilon$.
If $|\xi|+|\zeta|<\rho'$  then (\ref{Uee}) holds for all $\ve_n<\rho'$.
Since there are only finitely many functions $U^{\ve_n}$ with $\ve_n\geq \rho'$,
by choosing $\rho\in \,]0, \rho']$ small enough, we can guarantee that 
(\ref{Uee}) holds whenever $|\xi|+|\zeta|<\rho$ and $U^\ve$ 
is one of these finitely many functions with $\ve_n>\rho'$.
\v
{\bf 6.} Using the Riesz-Kolmogorov compactness theorem, 
we obtain a sequence $\ve\to 0$
such that
$$(u^\ve ,w^\ve,z^\ve,p^\ve,q^\ve)(X,Y)~\to~(u,w,z,p,q)(X,Y) \qquad\qquad\hbox{for
a.e.}~~(X,Y)\in \Omega^+.$$
By Lipschitz continuity, for any given $M>0$ this implies:
\begi
\item For a.e.~$X\in \R$ one has the uniform convergence
$(u^\ve ,w^\ve,p^\ve)(X,Y)~\to~(u,w,p)(X,Y)$,
 for all $Y\in[\vp(X), M]$.

\item For a.e.~$Y\in \R$ one has the uniform convergence
$(z^\ve,q^\ve)(X,Y)~\to~(z,q)(X,Y)$,
for all $X\in[\vp^{-1}(Y), M]$.
\endi
\v
{\bf 7.} It remains to
show that the limit functions $(u,w,z,p,q)$ provide a solution to the 
set of equations (\ref{2.17})--(\ref{2.20}).  This is nontrivial, because
the right hand side of the equations (\ref{2.17})--(\ref{2.20}) is
discontinuous at $w=\pi$ or $z=\pi$.  Comparing the functions $\theta,\theta_\ve$
in (\ref{thdef}) and  (\ref{thetadef}), one should be aware that 
general it is not true that
$\theta_\ve\to\theta$ as $\ve\to 0$.   For example, this convergence fails
if $w^\ve=z^\ve =\pi-\ve$.

Our proof is based on the 
following a priori estimate.
For any given $\eta,\ve>0$, consider the set
$$A_{\eta,\ve}~\doteq~
\bigl\{ (X,Y)\in \Omega_M\,;~~~\max\{w^\ve(X,Y), z^\ve(X,Y)\}\geq\pi-\eta
\bigr\}\,.$$
We claim that
\bel{ico}
\dint_{A_{\eta,\ve}} \Big( |w^\ve_Y|+|p^\ve_Y| + |z^\ve_X|+ |q^\ve_X|\Big)\, dXdY~
\leq~C\,
\eta^{1/3},\eeq
for some constant $C$, uniformly valid
 on the region where $0<\ve\le \eta$.
The idea of the proof is quite simple: on the set where $|w-z|$ is small
the derivatives $w_Y$, $p_Y$, $z_X$, $q_X$ are all close to zero.
On the other hand, the set where $|w-z|$ is large has small measure.
To simplify notation, we here omit the superscript $^\ve$.
More precisely, 
 for any $\delta>0$ consider the sets
$$\bega{rl}S_\delta&\doteq~\bigl\{ (X,Y)\in A_{\eta,\ve}\,;~~~w\leq z-\delta\leq z
\leq \pi
\bigr\}\,,\cr\cr
S'_\delta&\doteq~\bigl\{ (X,Y)\in A_{\eta,\ve}\,;~~~z\leq w-\delta\leq w
\leq \pi
\bigr\}\,.\enda$$ 
Observe that, for suitable constants $0<c<C$, we have 
\begi
\item[(i)]~$z_X (X,Y)\geq~~c\cdot \delta^2$ for all $(X,Y)\in S_\delta\,$.
\item[(ii)] If $z(X,Y)\geq \pi-\eta$, then 
$$
\left\{\bega{rl} z_X(X',Y)&\geq~-C\cdot \eta\,,\cr
 z(X',Y)&\geq~\pi-C \eta\,,\enda\right.\qquad
 \hbox{for all} ~~(X',Y)\in \Omega_M\,, ~~X'\geq X. $$
\endi
Entirely similar estimates hold for $w, w_Y$.
{}From (i)-(ii) we deduce
\bel{mS}
\hbox{meas}(S_\delta\cup S'_\delta)~=~\O(1)\cdot \eta\,\delta^{-2}.\eeq
Choosing $\delta = \eta^{1/3}$ we obtain
\bel{ic2}\bega{l}\ds
\left( \dint_{S_\delta\cup S'_\delta} 
+\dint_{A_{\eta,\ve}\setminus (S_\delta\cup S'_\delta)} 
\right)\Big( |w^\ve_Y|+|p^\ve_Y| + |z^\ve_X|+ |q^\ve_X|\Big)\, dXdY
\cr\cr
\qquad =~\O(1)\cdot \eta\, \delta^{-2} +\O(1)\cdot\delta ~=~\O(1)\cdot \eta^{1/3}.
\enda\eeq
\v
{\bf 8.} Thanks to (\ref{ico}), by 
choosing a further subsequence $\ve_n\downarrow 0$ and setting 
$\eta_n = \ve_n^{1/3}$, we can assume that
\bel{ic3}
\dint_{A_{\eta_n,\ve_m}} 
\Big( |w^{\ve_m}_Y|+|p^{\ve_m}_Y| + |z^{\ve_m}_X|+ |q^{\ve_m}_X|\Big)\, dXdY~\leq~C\,
\ve_n^{1/3}~\leq~2^{-n}\eeq
 for every $m\geq n\geq 1$.
 
Calling 
$$A_0~\doteq~
\bigl\{ (X,Y)\in \Omega_M\,;~~~\max\{w(X,Y), z(X,Y)\}\geq\pi
\bigr\},
$$
for every $n\geq 1$ we have the estimate
\bel{da0}\bega{l}
\ds\dint_{A_0} 
\Big( |w_Y|+|p_Y| + |z_X|+ |q_X|\Big)\,dXdY~\leq~ \dint_{A_{\eta_n,\ve_m}} 
\Big( |w_Y|+|p_Y| + |z_X|+ |q_X|\Big)\, dXdY\cr\cr
\qquad \leq~\limsup_{m\to\infty}
\dint_{A_{\eta_{n-1},\ve_m}} 
\Big( |w^{\ve_m}_Y|+|p^{\ve_m}_Y| + |z^{\ve_m}_X|+ |q^{\ve_m}_X|\Big)\, dXdY~
\leq~2^{-n+1}.
\enda\eeq
Hence the left hand side of (\ref{da0}) is zero.
\v
{\bf 9.} To complete the proof, setting 
$$B_0~\doteq~
\bigl\{ (X,Y)\in \Omega_M\,;~~~\max\{w(X,Y), z(X,Y)\}<\pi
\bigr\},
$$
 we need to show that 
\bel{ib0}\dint_{B_0} \Lambda(X,Y)\,dXdY~=~0,\eeq
where $\Lambda$ accounts for the differences between the right and left hand sides
of (\ref{2.17})--(\ref{2.20}). More precisely:
\bel{Ldef}\bega{rl}
\Lambda&\doteq~\ds
\left|w_Y-{c'(u)\over 8c^2(u)}\,( \cos z - \cos w)\,q\right| +\left|
z_X -{c'(u)\over 8c^2(u)}\,( \cos w - \cos z)\,p\right|\cr\cr
&\ds\qquad +\left| p_Y-{c'(u)\over 8c^2(u)}\,
(\sin z-\sin w)\,pq\right| +
\left|q_X-{c'(u)\over 8c^2(u)}\,
(\sin w-\sin z)\,pq\right|+\left|u_Y-{\sin z\over 4c}\right|\,.
\enda
\eeq
Toward this goal, for any $\nu\geq 1$ call
$$B_\nu~\doteq~
\bigl\{ (X,Y)\in \Omega_M\,;~~~\max\{w(X,Y), z(X,Y)\}\leq\pi-2^{-\nu}
\bigr\}.
$$
We can choose a sequence $\ve_n\downarrow0$ such that
$$(u^{\ve_n} ,w^{\ve_n},z^{\ve_n},p^{\ve_n},q^{\ve_n})(X,Y)~\to~(u,w,z,p,q)(X,Y) \qquad\qquad\hbox{for
a.e.}~~(X,Y)\in B_\nu.$$
By Egoroff's theorem, for any $\delta>0$, there exists a subset $F\subseteq B_\nu$ with $\mbox{meas}(F)<\delta$ and 
$$(u^{\ve_n} ,w^{\ve_n},z^{\ve_n},p^{\ve_n},q^{\ve_n})~\to~(u,w,z,p,q) \qquad\qquad\hbox{uniformly for any}~~(X,Y)\in B_\nu\backslash F.$$
Then for any $(X,Y)\in B_{\nu}\backslash F$, there exists some $\nu_1>\nu\geq 1$ such that
$$\max\{w^{\ve_n}(X,Y), z^{\ve_n}(X,Y)\}\leq\pi-2^{-\nu_1}.$$
It implies that
$$B_{\nu}\backslash F~ \subseteq ~B_{\nu_1,\ve_n}~\doteq~
\bigl\{ (X,Y)\in \Omega_M\,;~~~\max\{w^{\ve_n}(X,Y), z^{\ve_n}(X,Y)\}\leq\pi-2^{-\nu_1}
\bigr\}.
$$
Let $\Lambda^{\ve_n}$ be the same form of $\Lambda$ in (\ref{Ldef}) with $(u,w,z,p,q)$ replaced by 
$(u^{\ve_n} ,w^{\ve_n},z^{\ve_n},p^{\ve_n},q^{\ve_n})$. Equations (\ref{2.17e})-(\ref{2.20e}) implies that $\Lambda^{\ve_n}=2\ve_n$ on $B_{\nu_1,\ve_n}$.
Thus we obtain
\bel{ibnn}\bega{l}
\ds\dint_{B_\nu} \Lambda(X,Y)\,dXdY
=~\dint_{B_\nu\backslash F} \Lambda(X,Y)\,dXdY+\dint_{F} \Lambda(X,Y)\,dXdY
\cr\cr
\qquad\ds=~ \lim\limits_{\ve_n\downarrow0}\dint_{B_\nu\backslash F} \Lambda^{\ve_n}(X,Y)\,dXdY+\dint_{F} \Lambda(X,Y)\,dXdY
\cr\cr
\qquad\ds\leq~\lim\limits_{\ve_n\downarrow0}\dint_{B_{\nu_1,\ve_n}} \Lambda^{\ve_n}(X,Y)\,dXdY+C\delta
\cr\cr
\qquad=~C\delta.
\enda
\eeq
Hence 
\bel{ibn}
\dint_{B_\nu} \Lambda(X,Y)\,dXdY
=~0,
\eeq
for every $\nu\geq 1$.
Letting $\nu\to \infty$ in (\ref{ibn}) and using Lebesgue's monotone convergence theorem, we conclude (\ref{ib0}).

\endproof

\section{Global existence of dissipative solutions}
\label{sec:4}
\setcounter{equation}{0}

 Going back to the original variables, we now prove that the limit function $u$ 
provides a dissipative solution
to the original wave equation (\ref{1.1}).   
The proof of Theorem~2 will be given in several steps.

\v
{\bf 1.}  As in \cite{BZ}, by setting $f=x$ and then $f=t$ in (\ref{2.14}), we obtain 
\begin{equation}\label{th66}
\left\{
\begin{array}{ccc}
\displaystyle x_X=\frac{(1+\cos w)p}{4}\,,\\
\displaystyle  x_Y=-\frac{(1+\cos z)q}{4}\,,
\end{array}
\right.
\qquad\qquad
\left\{
\begin{array}{ccc}
\displaystyle t_X=\frac{(1+\cos w)p}{4c}\,,\\
\displaystyle  t_Y=\frac{(1+\cos z)q}{4c}\,.
\end{array}
\right.
\end{equation} 
Conversely, whenever $\cos w\neq-1$ and $\cos z\neq-1$, one has
\begin{equation}\label{th7}
\left\{
\begin{array}{ccc}
\displaystyle X_x=\frac{2}{(1+\cos w)p}\,,\\
\displaystyle  Y_x=-\frac{2}{(1+\cos z)q}\,,
\end{array}
\right.
\qquad\qquad
\left\{
\begin{array}{ccc}
\displaystyle X_t=\frac{2c}{(1+\cos w)p}\,,\\
\displaystyle  Y_t=\frac{2c}{(1+\cos z)q}\,.
\end{array}
\right.
\end{equation} 
As in \cite{BZ},
we can recover $(t,x)$ by integrating either one of the equations for  $x$ and
for $t$ in (\ref{th66}).
Indeed,  from the equations (\ref{2.17})--(\ref{2.20}) it follows
that  $x_{XY}=x_{YX}$ and $t_{XY}=t_{YX}$ for a.e.~$X,Y$. 

\begin{figure}[ht]
\centerline{\hbox{\includegraphics[width=10cm]{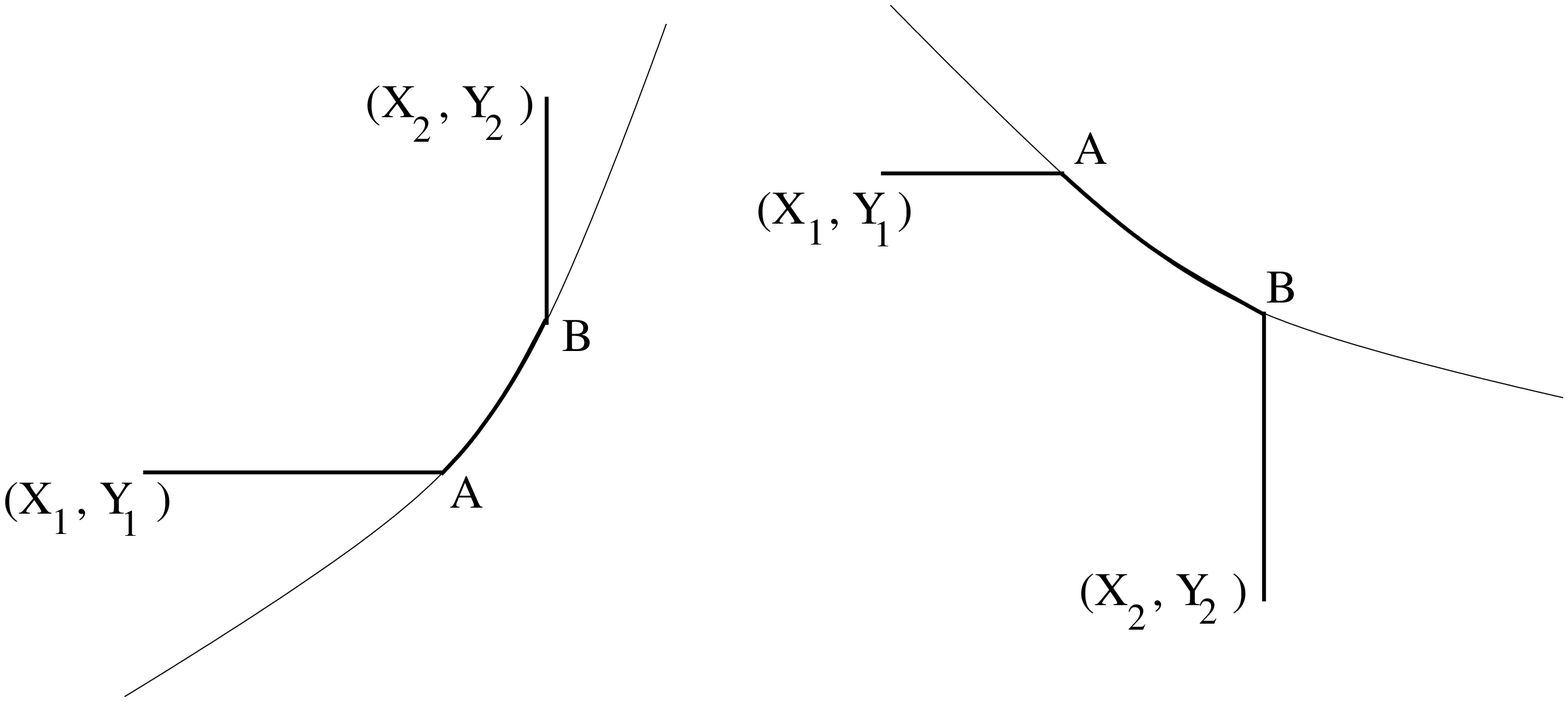}}}
\caption{\small Proving that the map $(t,x)\mapsto u(t,x)$ is well defined.}
\label{f:b84}
\end{figure}

\v
{\bf 2.}
For any $(\bar t,\bar x)$, we now define $u(\bar t,\bar x)\doteq u(X, Y)$ 
where $(X,Y)$ is any point such that $x(X,Y)=\bar x$ and $t(X,Y)=\bar t$. 
We  claim that the above definition of 
$u(t,x)$ is independent of the choice of $(X,Y)$. 
Indeed (see Fig.~\ref{f:b84}), suppose that there are two 
different points $(X_1,Y_1)$ and $(X_2,Y_2)$ such that
\begin{equation}\label{th8}
x(X_1,Y_1)~=~x(X_2,Y_2)~=~\bar x\, ,\qquad\qquad t(X_1,Y_1)~=~t(X_2,Y_2)~=~\bar t
\, . 
\end{equation}
Two cases must be considered.

CASE 1: $X_1\leq X_2$, $Y_1\leq Y_2$. We then consider the set
$$
\Gamma_{\bar x }~\doteq~\left\{(X,Y);~x(X,Y)\leq \bar x \right\}
$$
with boundary $\partial \Gamma_{\bar x }$. 
By (\ref{th66}), $x(X,Y)$ is increasing w.r.t.~$X$ and decreasing w.r.t.~$Y$.
This boundary can thus be represented as
the graph of a Lipschitz continuous function, namely $$X-Y~=~\phi(X+Y).$$
We now construct the Lipschitz continuous curve $\gamma$ as in Fig.~\ref{f:b84}, left, 
consisting
of
\begin{description}

\item[-] a horizontal segment joining $(X_1,Y_1)$ with a point $A=(X_A,Y_A)$
on $\partial \Gamma_{\bar x}$, with $Y_A=Y_1$,

\item[-] a portion of the boundary $\partial \Gamma_{\bar x}$,

\item[-] a vertical segment joining $(X_2,Y_2)$ to a point $B=(X_B,Y_B)$
on $\partial \Gamma_{\bar x}$, with $X_B=X_2$.
\end{description}
We can parameterize
this curve  in a Lipschitz continuous way, say $\gamma:[\xi_1,\xi_2]\mapsto \R^2$, using 
the parameter $\xi=X+Y$.
Observe that the map $(X,Y)\mapsto (t,x)$ is constant along $\gamma$.
By (\ref{th66}) this implies
$(1+\cos w) X_\xi=(1+\cos z)Y_\xi = 0$,
hence $\sin w\cdot X_\xi=\sin z\cdot Y_\xi=0$. We now compute
$$u(X_2,Y_2)-u(X_1,Y_1)~=~\int_\gamma \big(u_X\,dX+u_Y\,dY\big) ~=~\int_{\xi_1}^{\xi_2} \left({p\,\sin w\over 4c}\,X_\xi-{q\,\sin z\over 4c}\,
Y_\xi\right)\,d\xi~=~0\,,
$$
proving our claim.
\v
CASE 2: $X_1\leq X_2$, $Y_1\geq Y_2$.   In this case, we consider the
set
$$
\Gamma_{\bar t}\doteq \Big\{ (X,Y)\,;~~t(X,Y)\leq \bar t\Big\}\,,
$$
and construct a curve $\gamma$ connecting $(X_1,Y_1)$ with
$(X_2,Y_2)$ as in Fig.~\ref{f:b84}, right.  Details are entirely similar to Case 1.
\v
{\bf 3.} 
In this step
we prove that the function $u$ provides a weak solution
to the original nonlinear wave equation (\ref{1.1}). According to (\ref{1.5}),
we need to show that
\bel{ws1}
\dint \phi_t\big[(u_t+cu_x)+(u_t-cu_x)\big]
-\big(c(u)\phi\big)_x\big[(u_t+cu_x)-(u_t-cu_x)\big]
\,dxdt~=~0\eeq
for every test function $\phi\in \C^\infty_c([0,\infty[\,\times\R)$.
We now express the double integral in terms of the variables $X,Y$, using the 
change of variable formula
\bel{dxdt}
dxdt~=~{pq\over 2c(1+R^2)(1+S^2)}\, dXdY~=~{pq\over 2c} \cos^2{w\over 2} \cos^2{z\over 2}\, dXdY\,,
\eeq 
see \cite{BZ}
for details. Since only the absolutely continuous part of the measure
 $c'(u)u_x^2$ is accounted in the double integral (\ref{ws1}), using (\ref{2.14})
we obtain that (\ref{ws1}) is equivalent to 
 \bel{ws2}
\begin{split}
0~=&\dint R\big[\phi_t-(c\phi)_x)\big]
+S\big[\phi_t+(c\phi)_x\big]\,dxdt\\
~=&\dint -2cY_x\phi_YR+2cX_x\phi_XS+\theta c'\phi(u_XX_x+u_YY_x)(S-R)\,dxdt\,.
\end{split}\eeq
\v
{\bf Remark.} The integrand in (\ref{1.5}), or equivalently (\ref{ws2}), contains the 
term 
$c'(u)u_x^2$ which multiplies the test function $\phi$.   
Computing the same integral in terms of
$X,Y$, a straightforward use of the change of variable formula
would lead to a Radon measure.    However, it is only the absolutely continuous part
of this measure that actually contributes to the integral (\ref{1.5}), i.e., the 
part with $\max\{w,z\}<\pi$.
For this reason,  in (\ref{ws2}) we need to insert the additional factor $\theta$.
We observe that, in the conservative case \cite{BZ}, this factor was not needed, because
in that case the corresponding Radon measure is already absolutely continuous
with density $c'(u) u_x^2$, for a.e.~time $t$.
 \v
Using (\ref{dxdt}) to change variables and the identities
\bel{rswz}
\left\{
\bega{ccc}
\ds \frac{1}{1+R^2}~=~\cos^2\frac{w}{2}~=~\frac{1+\cos w}{2},\\
\\
\ds \frac{1}{1+S^2}~=~\cos^2\frac{z}{2}~=~\frac{1+\cos z}{2},
\enda
\right.\qquad\qquad 
\left\{
\bega{ccc}
\ds \frac{R}{1+R^2}~=~\frac{\sin w}{2},\\
\\
\ds \frac{S}{1+S^2}~=~\frac{\sin z}{2},
\enda
\right.
\eeq
the double integral in (\ref{ws2}) can be written as 
 \bel{ws3}
\begin{split}
&\dint \left\{\frac{R}{1+R^2}p\phi_Y+\frac{S}{1+S^2}q\phi_X+\theta\frac{c'pq}{8c^2}\left(\frac{\sin w}{1+S^2}-\frac{\sin z}{1+R^2}\right)(S-R)\phi\right\}\,dXdY\\
&=~\dint \left\{\frac{p\sin w}{2}\phi_Y+\frac{q\sin z}{2}\phi_X\right.\\
&\qquad\qquad \left.+\theta\frac{c'pq}{8c^2}\left(\sin w\sin z-\sin w\cos^2\frac{z}{2}\tan\frac{w}{2}
-\sin z\cos^2\frac{w}{2}\tan\frac{z}{2}\right)\phi\right\}\,dXdY\\
&=~\dint\left\{ {p\,\sin w\over 2}\,\phi_Y+{q\sin z\over 2}\,\phi_X +
\theta{c'pq\over 8c^2}\,\big[\cos (w-z)-1\big]\,\phi\right\} \,dXdY.
\end{split}\eeq
Since  $u,w,p$ are Lipschitz continuous functions of $Y$, while $u,z,q$ are 
Lipschitz continuous functions of $X$, after an integration by parts
it suffices to check that the identity
\bel{id1}
 \left({p\,\sin w\over 2}\right)_Y+\left({q\sin z\over 2}\right)_X ~=~
\theta{c'pq\over 8c^2}\,\big[\cos (w-z)-1\big]\eeq
holds at a.e.~point $(X,Y)$. By (\ref{2.17}) and (\ref{2.18}) we have
\bel{id2}
\begin{split}
 &\left({p\,\sin w\over 2}\right)_Y+\left({q\sin z\over 2}\right)_X\\ 
&=~p_Y{\sin w\over 2}+p\left({\sin w\over 2}\right)_Y+q_X{\sin z\over 2}+q\left({\sin z\over 2}\right)_X\\
&=~\theta\frac{c'pq}{16c^2}\big[(\sin z-\sin w)\sin w+\cos w(\cos z-\cos w)\big]\\
&\qquad\qquad +\theta\frac{c'pq}{16c^2}\big[(\sin w-\sin z)\sin z+\cos z(\cos w-\cos z)\big]\\
&=~\theta\frac{c'pq}{8c^2}\big[\cos(w-z)-1\big]
\end{split}
\eeq
which implies (\ref{id1}) and hence (\ref{ws1}). Therefore, the function $u$ provides a weak solution to (\ref{1.1}).

\begin{figure}[ht]
\centerline{\hbox{\includegraphics[width=10cm]{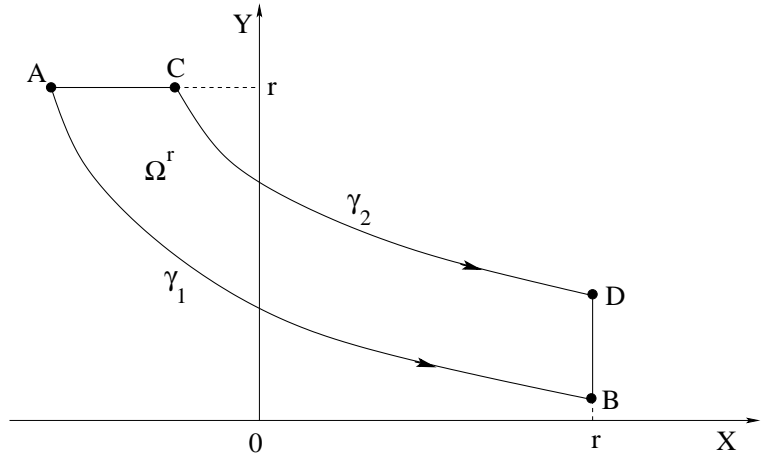}}}
\caption{\small The set $\Omega^r$ considered in step {\bf 4} of the proof.
In the contour integrations at (\ref{ch1}), 
the curve $\gamma_1$ is oriented from $A$ to $B$,
while $\gamma_2$ is oriented from $C$ to $D$.}
\label{f:wa19}
\end{figure}

\v
{\bf 4.}
It remains to prove that the weak solution $u$ is dissipative.
Toward this goal, consider any  $0\leq t_1<t_2$ and a large radius
$r>0$, and define
$$
\Omega^r~\doteq~\bigl\{(X,Y)\,;~~X\leq r,~~Y\leq r,~~ t_1\leq t(X,Y)\leq t_2\bigr\}.
$$
We can represent the above set as
$$\Omega^r~\doteq~\bigl\{(X,Y)\,;~~X\leq r,~~Y\leq r,~~ \phi_1(X+Y)\leq X-Y\leq
\phi_2(X+Y)\bigr\},
$$
for some functions $\phi_1<\phi_2$, Lipschitz continuous with constant 1.
With reference to Fig.~\ref{f:wa19}, assume that 
$$
x(A) = a, \qquad x(B)=b,\qquad x(C)=c,\qquad x(D)=d,$$
for some $a<b$ and $c<d$.
Moreover, let $\gamma_1$, $\gamma_2$ be the lower and upper portions of the boundary of
$\Omega^r$.
{}From the equations (\ref{2.17})-(\ref{2.18}) it follows that the 1-form 
\begin{equation}\label{th10}
E\,dx-(cM^2)\,dt~=~\frac{(1-\cos w)p}{8}\,dX-\frac{(1-\cos z)q}{8}\,dY
\end{equation}
is closed.
Therefore, the integral of the 1-form (\ref{th10}) along the boundary of $\Omega^r$ is 
zero. In particular, this yields
\bel{ch1}\bega{l}\ds
\int_{\gamma_1}\left\{\frac{(1-\cos w)p}{8}\,dX
-\frac{(1-\cos z)q}{8}\,dY\right\} - \int_{\gamma_2}\left\{\frac{(1-\cos w)p}{8}\,dX
-\frac{(1-\cos z)q}{8}\,dY\right\}\cr\cr
\qquad \ds=~\int_{AC}\frac{(1-\cos w)p}{8}\,dX
+\int_{BD}\frac{(1-\cos z)q}{8}\,dY~\geq~0\,.
\enda
\eeq

We now observe that, at time $t_1$, the total energy inside the interval $[a,b]$
is computed by
\bel{91}\bega{l}\ds
\int_a^b\frac{1}{2}\left[u^2_t(t_1,x)+c^2(u(t_1,x))u^2_x(t_1,x)\right]\,dx\cr\cr
\qquad\qquad \ds =~\int_{\gamma_1\cap\{ w(X,Y)<\pi\}}\frac{(1-\cos w)p}{8}\,dX
-\int_{\gamma_1\cap\{ z(X,Y)<\pi\}}\frac{(1-\cos z)q}{8}\,dY\,.
\enda\eeq
An entirely similar formula yields the energy at time $t_2$ inside the interval $[c,d]$.

\begin{figure}[ht]
\centerline{\hbox{\includegraphics[width=10cm]{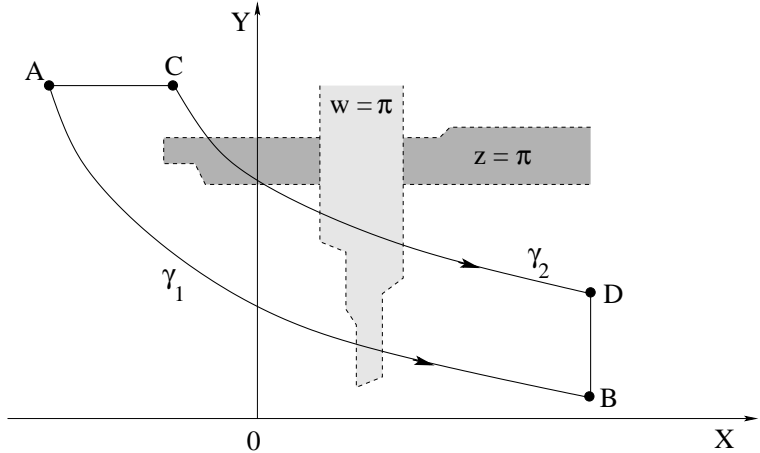}}}
\caption{\small  The set where $w=\pi$ grows as $Y$ increases.  Similarly, the set
where $z=\pi$ grows as $X$ increases.  
By (\ref{92}), this yields the inequalitites in (\ref{94}).}
\label{f:wa20}
\end{figure}

{}From the basic equations (\ref{2.17})--(\ref{thdef}) one obtains the implications
\bel{92}\bega{c}
w(X,Y)=\pi\quad\implies\quad w(X, Y')=\pi~~~
\hbox{and}~~~q(X,Y')=q(X,Y)\qquad\hbox{for all}~~Y'\geq Y,\cr
z(X,Y)=\pi\quad\implies\quad z(X', Y)~=~\pi~~~
\hbox{and}~~~p(X',Y)=p(X,Y)\qquad\hbox{for all}~~X'\geq X.\enda
\eeq
Combining (\ref{91}) with (\ref{ch1}) and using (\ref{92}) we now obtain
\bel{93}\bega{l}\ds\int_a^b\frac{1}{2}
\left(u^2_t(t_2,x)+c^2(u(t_2,x))u^2_x(t_2,x)\right)\,dx-
\int_c^d\frac{1}{2}\left(u^2_t(t_1,x)+c^2(u(t_1,x))u^2_x(t_1,x)\right)\,dx\cr\cr
\ds
~=~\left(\int_{\gamma_1}-\int_{\gamma_1\cap\{w(X,Y)=\pi\}}\right)
\frac{(1-\cos w)p}{8}\,dX-
\left(\int_{\gamma_1}-\int_{\gamma_1\cap\{z(X,Y)=\pi\}}\right)
\frac{(1-\cos z)q}{8}\,dY\cr\cr
\ds\qquad -\left(\int_{\gamma_2}-\int_{\gamma_2\cap\{w(X,Y)=\pi\}}\right)
\frac{(1-\cos w)p}{8}\,dX+
\left(\int_{\gamma_2}-\int_{\gamma_2\cap\{z(X,Y)=\pi\}}\right)
\frac{(1-\cos z)q}{8}\,dY\cr\cr
\ds ~=~
\int_{\gamma_1}\left\{\frac{(1-\cos w)p}{8}\,dX
-\frac{(1-\cos z)q}{8}\,dY\right\} - \int_{\gamma_2}\left\{\frac{(1-\cos w)p}{8}\,dX
-\frac{(1-\cos z)q}{8}\,dY\right\}\cr\cr
\quad \ds - \int_{\gamma_1\cap\{ w(X,Y)=\pi\}}{p\over 4}\, dX +
\int_{\gamma_1\cap\{ z(X,Y)=\pi\}}{q\over 4}\, dY
+ \int_{\gamma_2\cap\{ w(X,Y)=\pi\}}{p\over 4}\, dX -
\int_{\gamma_2\cap\{ z(X,Y)=\pi\}}{q\over 4}\, dY\cr\cr
~ \geq 0\,.
\enda\eeq
Indeed, as shown in Fig.~\ref{f:wa20}, by (\ref{92}) it follows
\bel{94}\bega{rl}0&\ds\leq~
\int_{\gamma_1\cap\{ w(X,Y)=\pi\}}{p\over 4}\, dX 
~\leq~\int_{\gamma_2\cap\{ w(X,Y)=\pi\}}{p\over 4}\, dX\,,\cr\cr
0&\ds \leq~- \int_{\gamma_1\cap\{ z(X,Y)=\pi\}}{q\over 4}\, dY~\leq~
-\int_{\gamma_2\cap\{ z(X,Y)=\pi\}}{q\over 4}\, dY.\enda\eeq

Letting $r\rightarrow+\infty$, we have $a,c\to -\infty$ while $b,d\to +\infty$.
Hence (\ref{93}) yields the desired inequality on the total energy:
$$\bega{l}\ds{\mathcal E}(t_2)~=~\int_{-\infty}^\infty\frac{1}{2}
\left(u^2_t(t_2,x)+c^2(u(t_2,x))u^2_x(t_2,x)\right)\,dx\cr\cr
\ds\qquad\leq~\int_{-\infty}^\infty\frac{1}{2}
\left(u^2_t(t_1,x)+c^2(u(t_1,x))u^2_x(t_1,x)\right)\,dx~=~
{\mathcal E}(t_1)\,,\enda$$ showing that the weak solution $u=u(t,x)$
is dissipative.
\endproof
\v
\begin{figure}[ht]
\centerline{\hbox{\includegraphics[width=13cm]{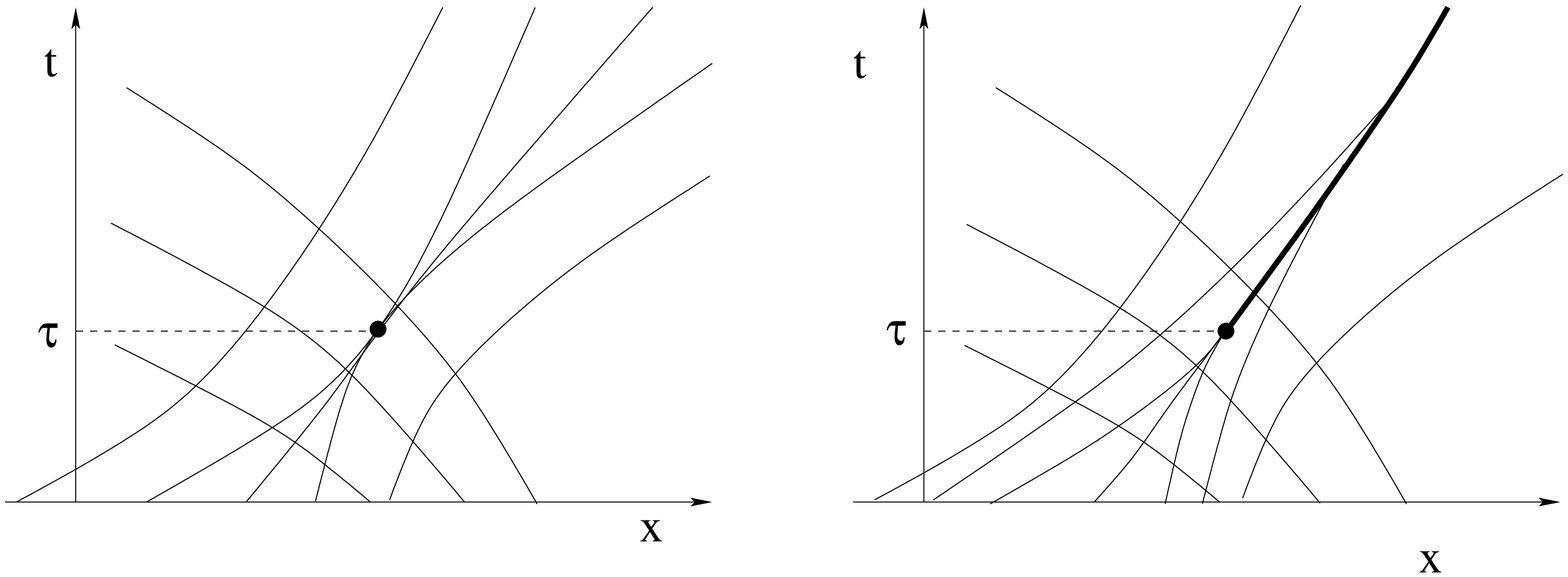}}}
\caption{\small Left: characteristic curves in a conservative solution.   At time $\tau$ a positive amount of energy is concentrated at a single point. However, for $t>\tau$
the energy measure is again absolutely  continuous.  Right: the characteristics curves
in a dissipative solution with the same initial data.  When some of 
the energy concentrates at one point, it remains inside the singular part of the 
energy measure $\mu^{(t)}$
for all subsequent times.    The equations (\ref{2.17})--(\ref{2.20})
imply that this singular part of the 
energy is formally transported along characteristics.   However
it does not affect the solution $u$ at any time  $t>\tau$.  }
\label{f:wa21}
\end{figure}
{\bf Remark.} Within the proof, we checked that the differential form 
(\ref{th10}) is closed.  This might suggest that, as  in 
\cite{BZ}, our solution is still conservative.   
The key difference can be explained as follows (see Fig.~\ref{f:wa21}).   
The contour integral 
$$\int_{\gamma_1}\left\{\frac{(1-\cos w)p}{8}\,dX
-\frac{(1-\cos z)q}{8}\,dY\right\} $$
yields the total energy at time $t_1$ inside the interval $[a,b]$.
In general, this energy is a positive Radon measure $\mu^{(t_1)}$ on the real line.   
On the other hand the integral
$$\int_a^b\frac{1}{2}\left[u^2_t(t_1,x)+c^2(u(t_1,x))u^2_x(t_1,x)\right]\,dx$$
accounts only for the absolutely continuous part of this measure.
In the solution constructed in \cite{BZ}, the measure $\mu^{(t)}$ is absolutely
continuous for a.e.~time $t$.  Hence ${\mathcal E}(t) = {\mathcal E(0)}$ for a.e.~$t$.
On the other hand, in our dissipative solution the singular part of the energy measure
(corresponding to $w=\pi$ or $z=\pi$)
is always increasing in time.   Hence the absolutely continuous part can only decrease.
\v
{\bf Acknowledgment.} This research was partially supported
by NSF, with grant  DMS-1411786: ``Hyperbolic Conservation Laws and Applications".

\end{document}